\documentclass[onecolumn,a4paper]{revtex4}
\usepackage[utf8]{inputenc}
\usepackage[english]{babel}

\usepackage{amsmath}    
\usepackage{verbatim}   
\usepackage{color}      
\usepackage{subfigure}  
\usepackage{hyperref}   
\usepackage{bbold}
\usepackage{mathrsfs}
\usepackage{textcomp}
\usepackage{amsthm}
\usepackage{amssymb}
\usepackage{mathtools}

\allowdisplaybreaks
\DeclareMathOperator\arccosh{arccosh}

\usepackage{tikz}
\usetikzlibrary{calc,decorations.markings}

\usepackage{scrextend}
\changefontsizes[15pt]{11pt}

 \usepackage[
    top    = 2.00cm,
    bottom = 2.00cm,
    left   = 1.75cm,
    right  = 1.75cm]{geometry}

\begin{document}
\title{Finite and infinite product transformations} 
\author{Martin Nicholson}
\begin{abstract}  
      Several infinite products are studied that satisfy the transformation relation of the type $f(\alpha)=f(1/\alpha)$. For certain values of the parameters these infinite products reduce to modular forms. Finite counterparts of these infinite products are motivated by solution of Dirichlet boundary problem on a rectangular grid. These finite product formulas give an elementary proof of several modular transformations.
\end{abstract}
\maketitle

$\mathbf{1}.$ B. Cais in an unpublished manuscript [\onlinecite{cais}] shows how to get the product transformation formula
\begin{equation}\label{modular1}
    \prod_{n=1}^\infty\left(\frac{1-e^{-\pi\alpha n}}{1+e^{-\pi\alpha n}}\right)^{(-1)^n}=\frac{1}{\sqrt{\alpha}}\prod_{n=1}^\infty\left(\frac{1-e^{-\pi n/\alpha}}{1+e^{-\pi n/\alpha}}\right)^{(-1)^n}
\end{equation}
from the well known series transformation [\onlinecite{berndt}]
\begin{equation}\label{sum}
    \sum_{n=1}^\infty \frac{n (-1)^n}{\sinh\pi\alpha n}+\alpha^{-2}\sum_{n=1}^\infty \frac{\alpha n(-1)^n}{\sinh\pi n/\alpha}=-\frac{1}{2\pi\alpha}
\end{equation}
by integration with respect to $\alpha$ using the indefinite integration
$$
\int\frac{ds}{\sinh s}=\ln\left(\frac{\sinh s}{1+\cosh s}\right).
$$
He also generalizes this approach to more complicated products like
\begin{equation}\label{modular2}
    \prod_{n=1}^\infty\left(1-\frac{2\sqrt{5}}{1+\sqrt{5}+4\cosh{\frac{2\pi\alpha n}{5}}}\right)^{\left(\frac{n}{5}\right)}=\prod_{n=1}^\infty\left(1-\frac{2\sqrt{5}}{1+\sqrt{5}+4\cosh{\frac{2\pi n}{5\alpha}}}\right)^{\left(\frac{n}{5}\right)},
\end{equation}
where $\left(\frac{n}{5}\right)$ denotes the Legendre symbol. Recall also the functional equation for the Dedekind eta product
\begin{equation}\label{dedekind}
e^{-\frac{\pi\alpha}{12}}\prod_{n=1}^\infty\left(1-e^{-2\pi\alpha n}\right)=\frac{1}{\sqrt{\alpha}}\cdot e^{-\frac{\pi}{12\alpha}}\prod_{n=1}^\infty\left(1-e^{-2\pi n/\alpha}\right).
\end{equation}
The aim of this paper is to generalize these transformation formulas and obtain finite products that reduce to the above formulas in the infinite limit. It should be noted that the infinite products in \eqref{modular1},\eqref{modular2},\eqref{dedekind} are modular forms, but the generalized products in the subsequent sections are not.

It is assumed throughout this paper that $\alpha>0$. However by analytic continuation all formulas are valid when $\text{Re}\,\alpha>0$.

$\mathbf{2}.$ Expanding $\text{cosech}$ into partial fractions and interchanging the order of summation we get
\begin{align}\label{trans1}
   f(\alpha,\gamma)&=\sum_{n=-\infty}^\infty \frac{\alpha(-1)^n}{\sqrt{\alpha^2n^2+\alpha\gamma^2}\sinh(\pi\sqrt{\alpha^2n^2+\alpha\gamma^2})}\\
    &\nonumber=\sum_{n=-\infty}^\infty{\pi\alpha(-1)^n}\sum_{m=-\infty}^\infty\frac{(-1)^m}{\pi^2(\alpha^2n^2+\alpha\gamma^2)+\pi^2 m^2}\\
    &\nonumber=\frac{1}{\pi}\sum_{m,n=-\infty}^\infty\frac{(-1)^{m+n}}{\alpha n^2+\alpha^{-1}m^2+\gamma^2}=f(1/\alpha,\gamma).
\end{align}
An identity equivalent to $f(\alpha,\gamma)=f(1/\alpha,\gamma)$ is given in \onlinecite{lattice} as equations (1.5.2) and (1.5.3) with $l=0$.

Now we multiply $ f(\alpha,\gamma)$ by $\pi\gamma$ and integrate termwise with respect to $\gamma$. Using the integral
\begin{align*}
&\int\frac{\pi\alpha\gamma\,d\gamma}{\sqrt{\alpha^2n^2+\alpha\gamma^2}\sinh(\pi\sqrt{\alpha^2n^2+\alpha\gamma^2})}\\
&=\int\frac{d(\pi\sqrt{\alpha^2n^2+\alpha\gamma^2})}{\sinh(\pi\sqrt{\alpha^2n^2+\alpha\gamma^2})}=\ln\left(\tanh\frac{\pi\sqrt{\alpha^2n^2+\alpha\gamma^2}}{2}\right)
\end{align*}
one can see that the sum
\begin{equation}\label{zerosum}
    \sum_{n=-\infty}^\infty (-1)^n\ln\left(\frac{\tanh({\pi\sqrt{\alpha^2n^2+\alpha\gamma^2}}/{2})}{\tanh({\pi\sqrt{\alpha^{-2}n^2+\alpha^{-1}\gamma^2}}/{2})}\right)
\end{equation}
doesn't depend on $\gamma$. Then it follows from the limit $\gamma\to\infty$  that \eqref{zerosum} is $0$. This means that
\begin{equation}\label{unitproduct}
    \prod_{n=-\infty}^\infty \left(\frac{\tanh({\pi\sqrt{\alpha^2n^2+\alpha\gamma^2}}/{2})}{\tanh({\pi\sqrt{\alpha^{-2}n^2+\alpha^{-1}\gamma^2}}/{2})}\right)^{(-1)^n}=1,
\end{equation}
or in the form easier to compare with \eqref{modular1}
\begin{equation}\label{nonmodular1}
    \prod_{n=1}^\infty\left(\frac{1-e^{-\pi\alpha\sqrt{n^2+\gamma^2}}}{1+e^{-\pi\alpha\sqrt{n^2+\gamma^2}}}\right)^{(-1)^n}=\sqrt{\frac{\tanh\frac{\pi\gamma}{2}}{\tanh\frac{\pi\alpha\gamma}{2}}}\prod_{n=1}^\infty\left(\frac{1-e^{-\pi\sqrt{n^2/\alpha^2+\gamma^2}}}{1+e^{-\pi\sqrt{n^2/\alpha^2+\gamma^2}}}\right)^{(-1)^n}.
\end{equation}
$\mathbf{3}.$ There is another way to write \eqref{trans1} in symmetric form. We start with
\begin{align*}
f(\alpha,\gamma)&=\sum_{n=-\infty}^\infty \frac{\alpha(-1)^n}{\sqrt{\alpha^2n^2+\alpha\gamma^2}}\cdot 2\sum_{m=0}^\infty e^{-\pi(2m+1)\sqrt{\alpha^2n^2+\alpha\gamma^2}}\\
&=2\sum_{m=0}^\infty \sum_{n=-\infty}^\infty\frac{(-1)^n}{\sqrt{n^2+\gamma^2/\alpha}} e^{-\pi(2m+1)\alpha\sqrt{n^2+\gamma^2/\alpha}}.
\end{align*}
and then apply Poisson summation formula to the sum over $n$ in the following form
\begin{align*}
\sum_{n=-\infty}^\infty\frac{(-1)^n}{\sqrt{n^2+y^2}} e^{-x\sqrt{n^2+y^2}}&=\sum_{n=-\infty}^\infty ~ \int\limits_{-\infty}^\infty \frac{e^{-\pi i t}}{\sqrt{t^2+y^2}} e^{-x\sqrt{t^2+y^2}}\cdot e^{-2\pi int}dt.
\end{align*}
The integral is calculated via formula $3.961.2$ from [\onlinecite{gr}] and equals $2K_0\left(y\sqrt{\pi^2(2n+1)^2+x^2}\right)$. Thus
\begin{align*}
f(\alpha,\gamma)&=4\sum_{m=0}^\infty\sum_{n=-\infty}^\infty  K_0\left(\frac{\gamma}{\sqrt{\alpha}}\sqrt{\pi^2(2n+1)^2+\pi^2(2m+1)^2\alpha^2}\right)\\
&=\sum_{m,n=-\infty}^\infty 2K_0\left(\pi\gamma\sqrt{\frac{(2n+1)^2}{\alpha}+(2m+1)^2\alpha}\right).
\end{align*}

$\mathbf{4.}$ The function
\begin{equation}
    h(z)=\frac{1}{\sqrt{z^2+\gamma^2}\left(1+e^{\pi\alpha \sqrt{z^2+\gamma^2}}\right)}\cdot\frac{1}{1+e^{-\pi iz}},
\end{equation}
where it is assumed that $\gamma>0$, is analytic in the complex plane with a cut $[i\gamma,+i\infty)$. It has simple poles on the real line at $z_n=2n+1$, and on the imaginary line at $\zeta_n=\sqrt{(2n+1)^2/\alpha^{2}+\gamma^2}$, $n\in\mathbb{Z}$. Consider contour $C$ in Fig.1, with small circles and semicircles of radii $\varepsilon$ around poles of $h(z)$. According to residue theorem $\int_Ch(z)dz=0$. Integrals along large arcs are $0$ in the limit $R\to \infty$. Sum of integrals along straight horizontal segments in the limit $R\to\infty$, $\varepsilon\to 0$ is
$$
I_h=2\int\limits_0^\infty \frac{dx}{\sqrt{x^2+\gamma^2}\left(1+e^{\pi\alpha \sqrt{x^2+\gamma^2}}\right)}.
$$
Similarly, sum of integrals along straight vertical segments in the limit $R\to\infty$, $\varepsilon\to 0$ is
$$
I_v=2\int\limits_{i\gamma}^{+i\infty} \frac{dz}{\sqrt{z^2+\gamma^2}\left(1+e^{-\pi iz}\right)}=-2\int\limits_{0}^{\infty} \frac{dy}{\sqrt{y^2+\gamma^2}\left(1+e^{\pi\sqrt{y^2+\gamma^2}}\right)}.
$$
Therefore in the limit $R\to\infty$, $\varepsilon\to 0$
$$
I_h+I_v-\pi i\sum_{n=-\infty}^\infty \underset{~z=z_n}{\text{res}}h(z)-2\pi i\sum_{n=0}^\infty \underset{~z=\zeta_n}{\text{res}}h(z)=0.
$$
\begin{figure}
    \centering
    \begin{tikzpicture}
    
  \draw [help lines,->] (-4,0) -- (4,0) coordinate (xaxis);
  \draw [help lines,->] (0,-0.5) -- (0,4) coordinate (yaxis);
  
  \draw [draw, line width=0.8pt] (0.8,0) arc (0:180:0.15);
  \draw [draw, line width=0.8pt] (-0.5,0) -- (0.5,0) ;
  \draw [draw, line width=0.8pt] (-0.5,0) arc (0:180:0.15);
  \draw [draw, line width=0.8pt] (0.8,0) -- (1.8,0) ;
  \draw [draw, line width=0.8pt] (2.1,0) arc (0:180:0.15);
  \draw [draw, line width=0.8pt] (2.1,0) -- (2.45,0) ;
  \draw [draw, line width=0.9pt] (2.75,0) node {$\ldots$};
  \draw [draw, line width=0.8pt] (3.05,0) -- (3.7,0) ;
  \draw [draw, line width=0.8pt] (3.7,0) arc (0:89:3.7);
  \draw [draw, line width=0.8pt] (0.061,3.7) -- (0.061,3.42) ;
  \draw [draw, line width=0.8pt] (-0.061,1.8) arc[radius = 0.15, start angle= 110, end angle= 250];
  \draw [draw, line width=0.8pt] (-0.061,1.52) -- (-0.061,0.99) ;
  \draw [draw, line width=0.8pt] (0.061,1.52) arc[radius = 0.15, start angle= 290, end angle= 430];
  \draw [draw, line width=0.8pt] (0.061,1.52) -- (0.061,0.99) ;
  \draw [draw, line width=0.8pt] (-0.061,1) arc[radius = 0.15, start angle= 112.5, end angle= 427.5];
  \draw [draw, line width=0.8pt] (-0.061,3.42) -- (-0.061,3.7) ;
  \draw [draw, line width=0.8pt] (-0.061,3.7) arc[radius = 3.7, start angle= 91, end angle= 180];
  \draw [draw, line width=0.8pt] (3.7,-0.25) node {$R$};
  \draw [draw, line width=0.8pt] (-3.7,-0.25) node {$-R$};
  \draw [draw, line width=0.8pt] (-3.05,0) -- (-3.7,0) ;
  \draw [draw, line width=0.9pt] (-2.75,0) node {$\ldots$};
  \draw [draw, line width=0.8pt] (-1.8,0) arc (0:180:0.15);
  \draw [draw, line width=0.8pt] (-1.8,0) -- (-0.8,0) ;
  \draw [draw, line width=0.8pt] (-2.45,0) -- (-2.1,0) ;
  \draw [draw, line width=0.9pt] (0.4,0.8) node {$i\gamma$};
  \draw [draw, line width=0.8pt] (0.061,1.8) -- (0.061,2.3) ;
  \draw [draw, line width=0.8pt] (-0.061,1.8) -- (-0.061,2.3) ;
  \draw [draw, line width=0.8pt] (0.061,2.3) arc[radius = 0.15, start angle= 290, end angle= 430];
  \draw [draw, line width=0.8pt] (-0.061,2.58) arc[radius = 0.15, start angle= 110, end angle= 250];
  \draw [draw, line width=0.8pt] (0.061,2.58) -- (0.061,2.96) ;
  \draw [draw, line width=0.8pt] (-0.061,2.58) -- (-0.061,2.96) ;
  \draw [draw, line width=0.9pt] (-0.061,3.3) node {$\vdots$};
  \draw [draw, line width=0.9pt] (0.061,3.3) node {$\vdots$};
  \draw [draw, line width=0.8pt] (0,-0.7) node {Fig.1};
  
\end{tikzpicture}
\end{figure}
Putting it altogether we arrive at
\begin{align*}
&\int\limits_0^\infty \frac{dx}{\sqrt{x^2+\gamma^2}\left(1+e^{\pi\alpha \sqrt{x^2+\gamma^2}}\right)}-\int\limits_{0}^{\infty} \frac{dy}{\sqrt{y^2+\gamma^2}\left(1+e^{\pi\sqrt{y^2+\gamma^2}}\right)}\\
&-\sum_{n=0}^\infty\frac{1}{\sqrt{(2n+1)^2+\gamma^2}\left(1+e^{\pi\alpha \sqrt{(2n+1)^2+\gamma^2}}\right)}+\sum_{n=0}^\infty\frac{1/\alpha}{\sqrt{(2n+1)^2/\alpha^2+\gamma^2}\left(1+e^{\pi\sqrt{(2n+1)^2/\alpha^2+\gamma^2}}\right)}=0.
\end{align*}
To convert it to a product one can multiply by $\pi\alpha\gamma$ and integrate with respect to $\gamma$ using
$$
\int\frac{\pi\alpha\gamma\,d\gamma}{\sqrt{x^2+\gamma^2}\left(1+e^{\pi\alpha \sqrt{x^2+\gamma^2}}\right)}=-\ln\left(1+e^{-\pi\alpha \sqrt{x^2+\gamma^2}}\right)
$$
to get the following symmetric relation
\begin{equation}\label{productintegral}
    \prod_{n=0}^\infty\frac{1+e^{-\pi\alpha\sqrt{(2n+1)^2+\gamma^2}}}{1+e^{-\pi\sqrt{(2n+1)^2/\alpha^2+\gamma^2}}}=\exp\left\{\frac{1}{2}\int\limits_0^\infty\ln\frac{1+e^{-\pi\alpha\sqrt{x^2+\gamma^2}}}{1+e^{-\pi\sqrt{x^2/\alpha^2+\gamma^2}}}\ dx\right\}.
\end{equation}

{\it{Note.}} Contour integrals similar to the one considered in this section have been investigated in ch. 1.5 of the book \onlinecite{lattice}, for example in eq. (1.5.1), with the aim of application to lattice sums. We study a more general integral in section 8.

$\mathbf{5}.$ The same analysis as in the preceding section is applicable to the function
\begin{equation}
    h(z)=\frac{1}{\sqrt{z^2+\gamma^2}\left(1-e^{\pi\alpha \sqrt{z^2+\gamma^2}}\right)}\cdot\frac{1}{1-e^{-\pi iz}}
\end{equation}
with the difference that the poles are now $z_n=2n$ and $\zeta_n=\sqrt{4n^2\alpha^{-2}+\gamma^2}$, $n\in\mathbb{Z}$. The result
\begin{equation}\label{dedekindgeneralization}
     \prod_{n=1}^\infty\frac{1-e^{-2\pi\alpha\sqrt{n^2+\gamma^2}}}{1-e^{-2\pi\sqrt{n^2/\alpha^2+\gamma^2}}}=\sqrt{\frac{1-e^{-2\pi\gamma}}{1-e^{-2\pi\alpha\gamma}}}\cdot \exp\left\{\int\limits_0^\infty\ln\frac{1-e^{-2\pi\alpha\sqrt{x^2+\gamma^2}}}{1-e^{-2\pi\sqrt{x^2/\alpha^2+\gamma^2}}}\ dx\right\},
\end{equation}
gives a function for which sum equals integral:
\begin{align}
\sum_{n=-\infty}^\infty\ln\frac{1-e^{-2\pi\alpha\sqrt{n^2+\gamma^2}}}{1-e^{-2\pi\sqrt{n^2/\alpha^2+\gamma^2}}}=\int\limits_{-\infty}^\infty\ln\frac{1-e^{-2\pi\alpha\sqrt{x^2+\gamma^2}}}{1-e^{-2\pi\sqrt{x^2/\alpha^2+\gamma^2}}}\, dx.
\end{align}
Analogous result is true for \eqref{productintegral}. In fact more generally
\begin{align}\label{jacobiimaginarygeneral}
    \nonumber\sum_{n=-\infty}^\infty\ln &\left\{\frac{1-e^{-2 \pi  \alpha  \sqrt{\gamma^2+(n+\theta)^2}}}{1-e^{-2 \pi  \sqrt{\gamma^2+{n^2}/{\alpha^2}}+2 \pi  i \theta}}\cdot\frac{1-e^{-2 \pi  \alpha  \sqrt{\gamma^2+(n-\theta)^2}}}{1-e^{-2 \pi  \sqrt{\gamma^2+{n^2}/{\alpha^2}}-2 \pi  i \theta}}\right\}\\
    &=\int\limits_{-\infty}^\infty \ln\left\{\frac{1-e^{-2 \pi  \alpha  \sqrt{\gamma^2+(x+\theta)^2}}}{1-e^{-2 \pi  \sqrt{\gamma^2+{x^2}/{\alpha^2}}+2 \pi  i \theta}}\cdot\frac{1-e^{-2 \pi  \alpha  \sqrt{\gamma^2+(x-\theta)^2}}}{1-e^{-2 \pi  \sqrt{\gamma^2+{x^2}/{\alpha^2}}-2 \pi  i \theta}}\right\}dx.
\end{align}
A list of functions for which sum equals integral along with references can be found in [\onlinecite{blog}].
Note that when $\gamma\to 0$ equation \eqref{dedekindgeneralization} reduces to the functional equation for the Dedekind eta function \eqref{dedekind}, while \eqref{jacobiimaginarygeneral} reduces to the Jacobi imaginary transform for theta functions [\onlinecite{ww}].

$\mathbf{6}.$ The choice
$$
h(z)=\frac{\cos\frac{3\pi z}{5}-\cos\frac{\pi z}{5}}{\sin \pi z}\cdot\frac{\cosh\frac{3\pi \alpha\sqrt{z^2+\gamma^2}}{5}-\cosh\frac{\pi \alpha\sqrt{z^2+\gamma^2}}{5}}{\sqrt{z^2+\gamma^2}\sinh (\pi \alpha \sqrt{z^2+\gamma^2})}
$$
with the help of formulas
$$
\cos\frac{3\pi n}{5}-\cos\frac{\pi n}{5}=(-1)^n\frac{\sqrt5}{2}\left(\frac{n}{5}\right),
$$
$$
\int\frac{\cosh\frac{3\pi t}{5}-\cosh\frac{\pi t}{5}}{\sinh \pi t}dt=\frac{\sqrt{5}}{2\pi}\ln\left(1-\frac{2\sqrt{5}}{1+\sqrt{5}+4\cosh{\frac{2\pi t}{5}}}\right)
$$
from [\onlinecite{cais}], leads to
\begin{equation}\label{legendre5}
    \prod_{n=1}^\infty\left(1-\frac{2\sqrt{5}}{1+\sqrt{5}+4\cosh{\frac{2\pi\alpha\sqrt{n^2+\gamma^2}}{5}}}\right)^{\left(\frac{n}{5}\right)}=\prod_{n=1}^\infty\left(1-\frac{2\sqrt{5}}{1+\sqrt{5}+4\cosh{\frac{2\pi\sqrt{n^2/\alpha^2+\gamma^2}}{5}}}\right)^{\left(\frac{n}{5}\right)}.
\end{equation}
This is generalization of \eqref{modular2}. From this it should be clear how to get generalized symmetric products from absolute invariants in section $5$ of [\onlinecite{cais}] automatically. However not all formulas seem to be amenable to such generalization (e.g.  see Proposition 25 in \cite{cais}).

$\mathbf{7}.$ Indeed, lets consider the function
$$
h(z)=\frac{\sin\frac{\pi  z}{3}  }{\sin \pi z}\cdot\frac{\sinh \frac{\pi \alpha\sqrt{z^2+\gamma^2}}{3}}{\sqrt{z^2+\gamma^2} \sinh\pi  \alpha\sqrt{z^2+\gamma^2}}.
$$
In this case, this function is integrated starting from the origin along first half of the contour depicted in Fig.1. and line segment $[0,i\gamma]$ connecting the point $i\gamma$ on the imagiary axis with the origin. It turns out that integrals along the real line and part of the imaginary line $(i\gamma,+i\infty)$ are real valued, while the integrals along small semicircles and the line segment $[0,i\gamma]$ are purely imaginary. Separating the imaginary part one arrives at the identity
\begin{align}\label{arctan1}
\nonumber\sum _{n=1}^{\infty} \left(\frac{n}{3}\right) \frac{\sinh \frac{\pi \alpha\sqrt{n^2+\gamma^2}}{3}}{\sqrt{n^2+\gamma^2} \sinh\pi  \alpha\sqrt{n^2+\gamma^2}}&+
\sum _{n=1}^{\infty} \left(\frac{n}{3}\right) \frac{\sinh \frac{\pi\sqrt{n^2/\alpha^2+\gamma^2}}{3}}{\sqrt{n^2+\alpha^2\gamma^2} \sinh\pi \sqrt{n^2/\alpha^2+\gamma^2}}\\&=
\frac{2}{\sqrt{3}}\int\limits_0^{\gamma} \frac{\sinh\frac{\pi t}{3}}{\sinh \pi  t} \cdot\frac{\sinh\frac{\pi \alpha  \sqrt{\gamma^2-t^2}}{3} }{\sqrt{\gamma^2-t^2}\,\sinh \pi  \alpha\sqrt{\gamma^2-t^2}}\, dt.
\end{align}
Now following the logic of the preceding sections we multiply this formula by $\gamma$ and integrate wrt $\gamma$. It is possible to do this integration in closed form for each term of the series on the lhs using the integration formula \cite{cais}
\begin{align*}
    \int_t^\infty \frac{\sinh s}{\sinh 3s} \, ds&=\frac{1}{\sqrt{3}} \left(\arctan\frac{1}{\sqrt{3}}-\arctan\frac{\tanh t}{\sqrt{3}}\right)\\
    &=\frac{1}{\sqrt{3}}\arctan\frac{\sqrt{3}}{1+2e^{2t}}.
\end{align*}
On the rhs we make the change of variable $t\to\sqrt{t}$ and introduce the notation $\tau=\gamma^2$. Then (\ref{arctan1}) becomes
$$
\sum _{n=1}^{\infty} \left(\frac{n}{3}\right) \arctan\frac{\sqrt{3}}{1+2 e^{\frac{2\pi \alpha}{3} \sqrt{{n^2}+\gamma^2}}}+\sum _{n=1}^{\infty} \left(\frac{n}{3}\right) \arctan\frac{\sqrt{3}}{1+2 e^{\frac{2\pi}{3} \sqrt{{n^2}/{\alpha^2}+\gamma^2}}}=\frac{\pi\alpha}{6}\int\limits_{\tau}^\infty I(s)\,ds
$$
where
$$
I(s)=\int\limits_0^{s} dt\,\frac{\sinh\frac{\pi \sqrt{t}}{3}}{\sqrt{t}\sinh \pi  \sqrt{t}} \cdot\frac{\sinh\frac{\pi \alpha  \sqrt{s-t}}{3} }{\sqrt{s-t}\,\sinh \pi  \alpha\sqrt{s-t}}.
$$
$I(s)$ is a convolution. Thus, if we write
$$
f(t)=\frac{\sinh\frac{\pi \sqrt{t}}{3}}{\sqrt{t}\sinh \pi  \sqrt{t}} ,\qquad g(t)=\frac{\sinh\frac{\pi \alpha\sqrt{t}}{3}}{\sqrt{t}\sinh \pi  \alpha\sqrt{t}},
$$
Then
\begin{align*}
\int\limits_{\tau}^\infty I(s)\,ds&=\int\limits_{0}^\infty I(s)\,ds-\int\limits_0^{\tau} I(s)\,ds\\&=
\int\limits_0^\infty ds\int\limits_0^s f(t)g(s-t)\,dt+\int\limits_0^\tau ds\int\limits_0^s f(t)g(s-t)\,dt\\
&=\int\limits_0^\infty f(t)\,dt\int\limits_t^\infty g(s-t)\,ds+\int\limits_0^\tau f(t)\,dt\int\limits_t^\tau g(s-t)\,ds
\\
&=\int\limits_0^\infty f(t)\,dt\cdot\int\limits_0^\infty g(s)\,ds+\int\limits_0^\tau f(t)\,dt\int\limits_0^{\tau-t} g(s)\,ds\\
&=\frac{1}{3\alpha}-\frac{2\sqrt{3}}{\alpha}\int\limits_0^\tau f(t)\arctan\frac{\tanh \frac{\pi\alpha\sqrt{\tau-t}}{3}}{\sqrt{3}}\,dt.
\end{align*}

Redefining $\alpha$ according to $\alpha\to \frac{3}{2\pi}\alpha$, and after this redefining $\gamma$ as $\gamma\to \gamma/\sqrt{\alpha}$ and introducing a new parameter $\beta=\frac{4\pi^2}{9\alpha}$ the identity takes the following symmetric form:

If $\alpha\beta=\frac{4\pi^2}{9}$, then
\begin{align}
\nonumber\sum _{n=1}^{\infty} \left(\frac{n}{3}\right) &\arctan\frac{\sqrt{3}}{1+2 e^{\sqrt{\alpha^2n^2+\alpha\gamma^2}}}+\sum _{n=1}^{\infty} \left(\frac{n}{3}\right) \arctan\frac{\sqrt{3}}{1+2 e^{\sqrt{\beta^2n^2+\beta\gamma^2}}}\\&=\frac{1}{3}\arctan\frac{\sqrt{3}}{1+2 e^{\gamma\sqrt{\beta}}}+\frac{2\gamma}{\sqrt{3\alpha}}\int\limits_0^1 \!\arctan\!\frac{\sqrt{3}}{1+2e^{\gamma\sqrt{\alpha(1-t^2)}}}\frac{dt}{1+2\cosh(\gamma t\sqrt{\beta})}.
\end{align}

Similarly, if $\alpha\beta=\frac{\pi^2}{4}$, then 
\begin{align}
\nonumber\sum _{n=1}^{\infty} \chi_4(n)\arctan e^{-\sqrt{\alpha^2n^2+\alpha\gamma^2}}+\sum _{n=1}^{\infty} &\chi_4(n)\arctan e^{-\sqrt{\beta^2n^2+\beta\gamma^2}}\\&=\frac{1}{2}\arctan e^{-\gamma\sqrt{\beta}}+\frac{\gamma}{4\sqrt{\alpha}}\int\limits_0^1 \!\frac{\arctan e^{-\gamma\sqrt{\alpha(1-t)}}}{\sqrt{t\,}\cosh(\gamma \sqrt{\beta t\,})}\,dt,
\end{align}
where $ \chi_4(n)=\sin\frac{\pi n}{2}$ is Dirichlet character modulo $4$.

$\mathbf{8}.$ Note that in section $6$ there were no integrals analogous to integrals $I_v,I_h$ from section $4$, thanks to the oddness of $h(z)$  they canceled out. With suitable function $h(z)$ this proof can be adapted to prove \eqref{nonmodular1} as well. In fact more is true:
\begin{equation}
     f(\alpha)=\prod _{n=1}^\infty \left(\frac{\cosh \left(\pi  \cos \theta \sqrt{n^2 \alpha^2+\alpha\gamma^2}\right)-\cos \left(\pi  n \alpha \sin \theta\right)}{\cosh \left(\pi  \cos \theta \sqrt{n^2 \alpha^2+\alpha\gamma^2}\right)+\cos \left(\pi  n \alpha \sin \theta\right)}\right)^{(-1)^n}=f\left(\frac{1}{\alpha}\right)\frac{\tanh\frac{\pi \gamma\cos\theta}{2\sqrt{\alpha}}}{\tanh\frac{\pi \gamma\sqrt{\alpha}\cos\theta}{2}}.
\end{equation}
Without presenting all the details of the proof we note the essential steps only. Take
\begin{equation*}
    h(z)=\frac{1}{\sqrt{z^2+\gamma^2}}\cdot\frac{\sinh\left({\pi\alpha}\cos\theta \sqrt{z^2+\gamma^2}\right)\cos\left({\pi\alpha}z\sin\theta\right)}{\cosh \left(2\pi \alpha \cos\theta \sqrt{z^2+\gamma^2}\right)-\cos(2\pi\alpha z\sin\theta)}\cdot \frac{1}{\sin\pi z}.
\end{equation*}
Clearly $h(z)$ is odd, hence there will not be any integral contributions analogous to $I_h$ and $I_v$ from section $3$. Besides $z_n=n$ this function has simple poles at
$$
\zeta_n=\frac{n\sin\theta}{\alpha}+i\cos\theta\sqrt{n^2/\alpha^2+\gamma^2},\quad n\in\mathbb{Z}.
$$
$\zeta_n$ satisfies the relation
$$
\zeta_n\sin\theta-n/\alpha=i\cos\theta\sqrt{\zeta_n^2+\gamma^2}.
$$
The poles $z_n$ are outside the contour of integration, while the poles $\zeta_n$ are inside.
Residues at $z_n=n$ are easily calculated and the sum over these residues can be converted to a product with the help of the integral
$$
\pi\alpha\cos\theta\int\frac{\sinh\left({\pi\alpha t}\cos\theta\right)\cos\left({\pi\alpha}\sin\theta\right)}{\cosh \left(2\pi \alpha t\cos\theta\right)-\cos(2\pi\alpha z\sin\theta)} dt=\frac14\ln\frac{\cosh\left({\pi\alpha t}\cos\theta\right)-\cos\left({\pi\alpha}\sin\theta\right)}{\cosh\left({\pi\alpha t}\cos\theta\right)+\cos\left({\pi\alpha}\sin\theta\right)}.
$$
Now we calculate the residues at $\zeta_n$:
\begin{align*}
  \underset{~z=\zeta_n}{\text{res}}h(z)&=\frac{1}{2\pi\alpha\sin\pi \zeta_n}\cdot\frac{\sinh\left({\pi\alpha}\cos\theta \sqrt{\zeta_n^2+\gamma^2}\right)\cos\left({\pi\alpha\zeta_n}\sin\theta\right)}{\zeta_n\cos\theta\sinh \left(2\pi \alpha \cos\theta \sqrt{\zeta_n^2+\gamma^2}\right)+\sqrt{\zeta_n^2+\gamma^2}\sin\theta\sin(2\pi\alpha \zeta_n\sin\theta)}\\&= \frac{1}{4\pi\alpha\sin\pi \zeta_n}\cdot\frac{(-1)^n}{\zeta_n\cos\theta+i\sin\theta\sqrt{\zeta_n^2+\gamma^2}}\\
  &=\frac{1}{4\pi\alpha\sin\pi \zeta_n}\cdot\frac{(-1)^{n-1}i}{\sqrt{n^2/\alpha^2+\gamma^2}}. 
\end{align*}
From this it follows that
$$
\underset{~z=\zeta_n}{\text{res}}h(z)+\underset{~z=\zeta_{-n}}{\text{res}}h(z)=\frac{(-1)^{n-1}}{\pi\alpha\sqrt{n^2/\alpha^2+\gamma^2}}\cdot\frac{\sinh\left({\pi}\cos\theta \sqrt{n^2/\alpha^2+\gamma^2}\right)\cos\left({\pi}z\alpha^{-1}\sin\theta\right)}{\cosh \left(2\pi\cos\theta \sqrt{n^2/\alpha^2+\gamma^2}\right)-\cos(2\pi z\alpha^{-1}\sin\theta)}.
$$
The RHS of this expression has the same form as $\underset{~z=z_n}{\text{res}}h(z)$, and therefore can be converted to a product by integration.

$\mathbf{9}.$ The result of this section is quite similar to the previous section:
\begin{equation}\label{product_cosh}
    f(\alpha)=\prod _{n=0}^\infty \left(\frac{\cosh \left(\pi  \cos \theta \sqrt{\left(n+\frac{1}{2}\right)^2 \alpha^2+\alpha \gamma^2}\right)-\sin \left(\pi \alpha \left(n+\frac{1}{2}\right) \sin \theta\right)}{\cosh \left(\pi  \cos \theta \sqrt{\left(n+\frac{1}{2}\right)^2 \alpha^2+\alpha \gamma^2}\right)+\sin \left(\pi \alpha \left(n+\frac{1}{2}\right) \sin \theta\right)}\right)^{(-1)^n}=f\left(\frac{1}{\alpha}\right).
\end{equation}
In this case
\begin{equation*}
    h(z)=\frac{1}{\sqrt{z^2+\gamma^2}}\cdot\frac{\sinh\left({\pi\alpha}\cos\theta \sqrt{z^2+\gamma^2}\right)\sin\left({\pi\alpha}z\sin\theta\right)}{\cosh \left(2\pi \alpha \cos\theta \sqrt{z^2+\gamma^2}\right)+\cos(2\pi\alpha z\sin\theta)}\cdot \frac{1}{\cos\pi z},
\end{equation*}
and
$$
\pi\alpha\cos\theta\int\frac{\sinh\left({\pi\alpha t}\cos\theta\right)\sin\left({\pi\alpha}\sin\theta\right)}{\cosh \left(2\pi \alpha t\cos\theta\right)+\cos(2\pi\alpha z\sin\theta)} dt=\frac14\ln\frac{\cosh\left({\pi\alpha t}\cos\theta\right)-\sin\left({\pi\alpha}\sin\theta\right)}{\cosh\left({\pi\alpha t}\cos\theta\right)+\sin\left({\pi\alpha}\sin\theta\right)}.
$$

Note that $f\left(\frac{2N}{\sin\theta}\right)=1$ trivially for $N\in\mathbb{N}$. However $f\left(\frac{\sin\theta}{2N}\right)=1$ is quite non-trivial. As an example, by letting $N=1$, $\gamma=0$ in (\ref{product_cosh}) one obtains the identity
\begin{equation*}
\prod _{n=1}^\infty \left(\frac{\cosh\frac{\pi n\sin\theta\cos\theta}{4} -\sin\frac{\pi n\sin^2\theta}{4}}{\cosh\frac{\pi n\sin\theta\cos\theta}{4} +\sin\frac{\pi n\sin^2\theta}{4}}\right)^{\chi_4(n)}=1,
\end{equation*}
where $\chi_4(n)=\sin\frac{\pi n}{2}$ is Dirichlet character modulo $4$.

$\mathbf{10}.$ So far all products or series have been symmetric. But there are also identities for non-symmetric products as well. Here is one such example without proof:

\begin{equation}
    \prod_{n=-\infty}^\infty\frac{\tanh\pi\sqrt{\alpha^2n^2+\tfrac14}}{\Big(1-e^{-\tfrac{\pi}{\alpha}\sqrt{n^2+1}}\Big)^{(-1)^n}}=\exp\left\{\frac{1}{\alpha}\int\limits_{-\infty}^\infty\ln\left(\tanh\pi\sqrt{t^2+\tfrac14}\right)dt\right\}.
\end{equation}
This can be reformulated as the expression
$$
\left\{\prod_{n=-\infty}^\infty\frac{\tanh\pi\sqrt{\alpha^2n^2+\tfrac14}}{\Big(1-e^{-\tfrac{\pi}{\alpha}\sqrt{n^2+1}}\Big)^{(-1)^n}}\right\}^{\alpha}
$$
being independent of $\alpha$.

$\mathbf{11}.$ The coupled equations defining $\alpha_j$ and $\beta_k$ 
\begin{equation}\label{coupled1}
    \cosh\alpha_j+\cos\frac{\pi (j-\frac{1}{2})}{n}= \cosh\beta_k+\cos\frac{\pi (k-\frac{1}{2})}{m}\qquad(1\le j\le n,~1\le k\le m)
\end{equation}
arise in the solution of Laplace and Helmholtz equations on a lattice [\onlinecite{wiener}]. These set of $\alpha_j$ and $\beta_k$ satisfy the reciprocal relation
\begin{equation}\label{reciprocal1}
    \prod_{j=1}^n2\cosh m\alpha_j=\prod_{k=1}^m2\cosh n\beta_k.
\end{equation}
The proof is very simple and uses the well known formula (see [\onlinecite{jolley}], formula $1026$):
\begin{equation}\label{cosproduct}
    2^{m-1} \prod _{j=1}^m \left(\cosh \alpha-\cos \frac{\pi  (j-1/2)}{m}\right)=\cosh m \alpha.
\end{equation}
Indeed, denoting by $x$ the common value of the equations \eqref{coupled1}
\begin{align*}
    \prod_{j=1}^n2\cosh m\alpha_j&=\prod_{j=1}^n2\cdot 2^{m-1} \prod _{k=1}^m \left(\cosh \alpha_j-\cos \frac{\pi  (k-1/2)}{m}\right)\\
    &=2^{mn}\prod_{j=1}^n \prod _{k=1}^m  \left(x-\cos \frac{\pi  (j-1/2)}{n}-\cos \frac{\pi  (k-1/2)}{m}\right).
\end{align*}
This expression is symmetric in $m$ and $n$ and therefore imply \eqref{reciprocal1}.

$\mathbf{12}.$ Similarly, if
\begin{equation}\label{couple2}
   \cosh\alpha_j+ \cos\frac{\pi j}{n}= \cosh\beta_k+\cos\frac{\pi k}{m}\qquad(1\le j\le n,~1\le k\le m)
\end{equation}
then
\begin{equation}\label{reciprocal2}
    \prod_{j=1}^n\frac{\sinh m\alpha_j}{\sinh\alpha_j}=\prod_{k=1}^m\frac{\sinh n\beta_k}{\sinh\beta_k}.
\end{equation}
The proof is similar to the one in previous section, this time using the well known formula
\begin{equation}\label{sinproduct}
    2^{m-1} \prod _{k=1}^{m-1} \left(\cosh \alpha-\cos \frac{\pi  k}{m}\right)=\frac{\sinh m \alpha}{\sinh \alpha}.
\end{equation}

$\mathbf{13}.$ If
\begin{equation}
    \cosh\frac{\alpha_j}{2}\cos\frac{\pi (j-1/2)}{2n}=\cosh\frac{\beta_k}{2}\cos\frac{\pi (k-1/2)}{2m}\qquad(1\le j\le n,~1\le k\le m)
\end{equation}
then
\begin{equation}
    \prod_{j=1}^n\cosh m\alpha_j=\prod_{k=1}^m\cosh n\beta_k.
\end{equation}
In this case, besides \eqref{cosproduct} one needs 
\begin{equation}
    \prod_{j=1}^n \cos\frac{\pi (j-1/2)}{2n}=\frac{\sqrt2}{2^n}.
\end{equation}

$\mathbf{14}.$ Here is a less trivial example inspired by the `dispersion relation' $~\cosh \frac{\omega}{2}=\frac{2}{\cos p}-\cos p~$ on the triangular and hexagonal lattices [\onlinecite{cserti}]:
If
\begin{align}\label{coupled4}
    \left(\cosh\frac{\alpha_j}{2}+\cos\frac{\pi (2j-1)}{4n}\right)\cos\frac{\pi (2j-1)}{4n}=\left(\cosh\frac{\beta_k}{2}+\cos\frac{\pi (2k-1)}{4m}\right){\cos\frac{\pi (2k-1)}{4m}}
\end{align}
for $1\le j\le n,~1\le k\le m$, then
\begin{equation}\label{reciprocal4}
     \prod_{j=1}^n\left(\cosh m\alpha_j+\cos\frac{m\pi (2j-1)}{2n}\right)=\prod_{k=1}^m\left(\cosh n\beta_k+\cos\frac{n\pi (2k-1)}{2m}\right).
\end{equation}

To prove this relation, according to the following generalization of \eqref{cosproduct}
\begin{equation}
    2^{m-1} \prod _{j=1}^m \left[\cosh \alpha-\cos \left(y+\frac{2\pi j}{m}\right)\right]=\cosh m \alpha-\cos my,
\end{equation}
we write
\begin{align*}
   &\prod_{j=1}^n\left(\cosh m\alpha_j+\cos\frac{m\pi (2j-1)}{2n}\right)=    \prod_{j=1}^n2^{m-1} \prod _{k=1}^m\left[\cosh \alpha_j-\cos\left(\frac{\pi (2j-1)}{2n}+\frac{\pi (2k-1)}{m}\right)\right]\\
   & = 2^{n(2m-1)}\prod _{k=1}^m\prod_{j=1}^n \left[\cosh\frac{\alpha_j}{2}-\cos\left(\frac{\pi (2j-1)}{4n}+\frac{\pi (2k-1)}{2m}\right)\right]\left[\cosh \frac{\alpha_j}{2}+\cos\left(\frac{\pi (2j-1)}{4n}+\frac{\pi (2k-1)}{2m}\right)\right].
\end{align*}
Denote the common value of equations (\ref{coupled4}) by $x$. Then
\begin{align*}
     \prod_{j=1}^n&\left[\cosh\frac{\alpha_j}{2}-\cos\left(\frac{\pi (2j-1)}{4n}+\frac{\pi (2k-1)}{2m}\right)\right]\\
     &\phantom{...............}=\prod_{j=1}^n\left[\frac{x}{\cos\frac{\pi (2j-1)}{4n}}-\cos\frac{\pi (2j-1)}{4n}-\cos\left(\frac{\pi (2j-1)}{4n}+\frac{\pi (2k-1)}{2m}\right)\right]\\
    &\phantom{...............}=\frac{2^n}{\sqrt2}\prod_{j=1}^n\left[x-2\cos\frac{\pi (2j-1)}{4n}\cos\frac{\pi (2k-1)}{4m}\cos\left(\frac{\pi (2j-1)}{4n}+\frac{\pi (2k-1)}{4m}\right)\right],
\end{align*}
\begin{align*}
    \prod_{j=1}^n&\left[\cosh\frac{\alpha_j}{2}+\cos\left(\frac{\pi (2j-1)}{4n}+\frac{\pi (2k-1)}{2m}\right)\right]\\&\phantom{...............}=\prod_{j=1}^n\left[\frac{x}{\sin\frac{\pi (2j-1)}{4n}}-\sin\frac{\pi (2j-1)}{4n}+\sin\left(\frac{\pi (2j-1)}{4n}-\frac{\pi (2k-1)}{2m}\right)\right]\\
&\phantom{...............}=\frac{2^n}{\sqrt2}\prod_{j=1}^n\left[x-2\sin\frac{\pi (2j-1)}{4n}\sin\frac{\pi (2k-1)}{4m}\cos\left(\frac{\pi (2j-1)}{4n}-\frac{\pi (2k-1)}{4m}\right)\right].
\end{align*}
From these two formulas, it is obvious that the product on the LHS of \eqref{reciprocal4} is symmetric in $n$ and $m$. This completes the proof.

$\mathbf{15}.$
If 
\begin{equation}\label{coupled5}
    \cosh\alpha_j+\cos\frac{\pi j}{2n}=\cosh\beta_k+\cos\frac{\pi k}{2m},\qquad (1\le j\le 2n,\ 1\le k\le 2m),
\end{equation}
then
\begin{equation}\label{reciprocal5}
    \prod_{j=1}^{2n}\left(\frac{\tanh m\alpha_j}{\sinh\alpha_j}\right)^{(-1)^j}=\prod_{k=1}^{2m}\left(\frac{\tanh n\beta_k}{\sinh\beta_k}\right)^{(-1)^k}.
\end{equation}

Denote the common value of equations \eqref{coupled5} by $x$. First, note that $\alpha_{2n}=\beta_{2m}$ and as a result (\ref{reciprocal5}) is equivalent to
\begin{equation*}
    \prod_{j=1}^{2n-1}\left(\frac{\tanh m\alpha_j}{\sinh\alpha_j}\right)^{(-1)^j}=\prod_{k=1}^{2m-1}\left(\frac{\tanh n\beta_k}{\sinh\beta_k}\right)^{(-1)^k}.
\end{equation*}
Here we use \eqref{cosproduct} and \eqref{sinproduct} to write the product on the lhs in symmetric form:
\begin{align*}
\prod_{j=1}^{2n-1}\left(\frac{\tanh m\alpha_j}{\sinh\alpha_j}\right)^{(-1)^j}&=\frac{\sinh\alpha_{2n-1}}{\tanh m\alpha_{2n-1}}\prod_{j=1}^{n-1}\frac{\tanh m\alpha_{2j}}{\sinh\alpha_{2j}}\frac{\sinh\alpha_{2j-1}}{\tanh m\alpha_{2j-1}}\\
&=\frac{\sinh\alpha_{2n-1}}{\tanh m\alpha_{2n-1}}\prod_{j=1}^{n-1}\frac{\cosh\alpha_{2j-1}-\cos\frac{\pi (m-1/2)}{m}}{\cosh\alpha_{2j}-\cos\frac{\pi(m-1/2)}{m}}\\
&\times \prod_{j=1}^{n-1}\prod_{k=1}^{m-1}\frac{\cosh\alpha_{2j-1}-\cos\frac{\pi (k-1/2)}{m}}{\cosh\alpha_{2j}-\cos\frac{\pi(k-1/2)}{m}}\frac{\cosh\alpha_{2j}-\cos\frac{\pi k}{m}}{\cosh\alpha_{2j-1}-\cos\frac{\pi k}{m}}\\
&=\frac{1}{\cosh\alpha_{2n-1}-\cos\frac{\pi (m-1/2)}{m}}\cdot\frac{\sinh\alpha_{2n-1}}{\tanh m\alpha_{2n-1}}\frac{\sinh\beta_{2m-1}}{\tanh n\beta_{2m-1}}\\
&\times \prod_{j=1}^{n-1}\prod_{k=1}^{m-1}\frac{x-\cos\frac{\pi (j-1/2)}{n}-\cos\frac{\pi (k-1/2)}{m}}{x-\cos\frac{\pi j}{n}-\cos\frac{\pi(k-1/2)}{m}}\frac{x-\cos\frac{\pi j}{n}-\cos\frac{\pi k}{m}}{x-\cos\frac{\pi (j-1/2)}{n}-\cos\frac{\pi k}{m}}.
\end{align*}
Now, in this last formula the factor 
$$
\frac{1}{\cosh\alpha_{2n-1}-\cos\frac{\pi (m-1/2)}{m}}=\frac{1}{x+\cos\frac{\pi}{2n}+\cos\frac{\pi }{2m}}
$$
is obviously symmetric when $m$ and $n$ are interchanged, and so is the factor $\frac{\sinh\alpha_{2n-1}}{\tanh m\alpha_{2n-1}}\frac{\sinh\beta_{2m-1}}{\tanh n\beta_{2m-1}}$ and the double product. So both sides in $(4)$ are symmetric when $m$ and $n$ are interchanged, and hence they are equal.

$\mathbf{16}.$ What is the limit $n,m\to\infty$ of the reciprocal relation \eqref{reciprocal5}, if it exists? If the common value $x$ of equations \eqref{coupled5} is chosen to be close to $2$, then $\alpha_j$ will be close to $0$ for small $j$. Let $x=2+\frac{\pi^2\gamma^2}{8n^2}$, then expanding $\cos$ and $\cosh$ we get approximately for small $j$ and $k$
$$
\alpha_j=\frac{\pi\sqrt{ j^2+\gamma^2}}{2n},\quad \beta_k=\frac{1}{2n}\sqrt{\frac{n^2}{m^2}k^2+\gamma^2}.
$$
Next assume that $m,n\to\infty$ such that $\frac{m}{n}\to\alpha$. With these assumptions we have
\begin{equation}\label{limit1}
    \lim_{\substack{\scriptstyle{n,m\to\infty}\\m/n\to\alpha}} \prod_{j=1}^{2n-1}\left({\tanh m\alpha_j}\right)^{(-1)^j}=\prod_{j=1}^\infty\left(\tanh \frac{\pi\alpha\sqrt{ j^2+\gamma^2}}{2}\right)^{(-1)^j},
\end{equation}
\begin{equation}\label{limit2}
    \lim_{\substack{n,m\to\infty\\m/n\to\alpha}} \prod_{k=1}^{2m-1}\left({\tanh n\beta_k}\right)^{(-1)^k}=\prod_{k=1}^\infty\left(\tanh \frac{\pi\sqrt{ k^2/\alpha^2+\gamma^2}}{2}\right)^{(-1)^k}.
\end{equation}
The $\sinh$ factors are simplified as follows
\begin{align*}
    \prod_{j=1}^{2n-1}\left({\sinh \alpha_j}\right)^{(-1)^j}&=\prod_{j=1}^{2n-1}\left[(\cosh \alpha_j-1)(\cosh \alpha_j+1)\right]^{(-1)^j/2}\\
&=\prod_{j=1}^{2n-1}\left[\left(x-1-\cos\frac{\pi j}{2n}\right)\left(x+1-\cos\frac{\pi j}{2n}\right)\right]^{(-1)^j/2}\\
&=\left\{\frac{\tanh[n\arccosh(x-1)]}{\sinh[\arccosh(x-1)]}\frac{\tanh[n\arccosh(x+1)]}{\sinh[\arccosh(x+1)]}\right\}^{1/2}.
\end{align*}
Now it is easy to calculate the limit
\begin{align}\label{limit3}
    \nonumber&\lim_{\substack{n,m\to\infty\\m/n\to\alpha}} \frac{\prod_{j=1}^{2n-1}\left({\sinh \alpha_j}\right)^{(-1)^j}}{\prod_{k=1}^{2m-1}\left({\sinh \beta_k}\right)^{(-1)^k}}\\
    &=\lim_{\substack{n,m\to\infty\\m/n\to\alpha}} \left\{\frac{\tanh[n\arccosh(x-1)]\tanh[n\arccosh(x+1)]}{\tanh[m\arccosh(x-1)]\tanh[m\arccosh(x+1)]}\right\}^{1/2}=\sqrt{\frac{\tanh\frac{\pi\gamma}{2}}{\tanh\frac{\pi\alpha\gamma}{2}}}.
\end{align}
Combining equations \eqref{limit1}-\eqref{limit3} leads to \eqref{nonmodular1}. Thus we obtained an elementary proof of the modular transformation \eqref{modular1}. Elementary proofs of modular transformations are known in the literature. For example Apostol [\onlinecite{apostol}] gives an elementary proof for $\sum_{n=1}^\infty \frac{n^{-p}x^p}{1-x^p}$, $|x|<1$ for an odd integer $p>1$.

$\mathbf{17}.$
If 
\begin{equation}
    \cosh\alpha_j+\cos\frac{2\pi j}{5n}=\cosh\beta_k+\cos\frac{2\pi k}{5m},\qquad (1\le j\le 5n,\ 1\le k\le 5m),
\end{equation}
then
\begin{equation}
    \prod_{j=1}^{5n}\left(1-\frac{2\sqrt{5}}{1+\sqrt{5}+4\cosh m\alpha_j}\right)^{\left(\frac{j}{5}\right)}=\prod_{k=1}^{5m}\left(1-\frac{2\sqrt{5}}{1+\sqrt{5}+4\cosh n\beta_k}\right)^{\left(\frac{k}{5}\right)},
\end{equation}
or alternatively
\begin{equation}
    \prod_{j=1}^{5n}\frac{\cosh m\alpha_j-\cos\frac{2\pi j}{5}}{\cosh m\alpha_j-\cos\frac{4\pi j}{5}}=\prod_{k=1}^{5m}\frac{\cosh n\beta_k-\cos\frac{2\pi k}{5}}{\cosh n\beta_k-\cos\frac{4\pi k}{5}}.
\end{equation}

This is finite version of (\ref{legendre5}). Its proof is similar to the proof of equation (\ref{reciprocal5}) in section $14$.

$\mathbf{18}.$ If 
\begin{equation}
    \cosh\alpha_j+\cos\frac{\pi j}{6n}=\cosh\beta_k+\cos\frac{\pi k}{6m},\qquad (1\le j\le 6n,\ 1\le k\le 6m),
\end{equation}
then
\begin{equation}
    \prod_{j=1}^{6n}\left(\frac{2\cosh m\alpha_j-\sqrt{3}}{2\cosh m\alpha_j+\sqrt{3}}\right)^{\chi_{12}(j)}= \prod_{k=1}^{6m}\left(\frac{2\cosh n\beta_k-\sqrt{3}}{2\cosh n\beta_k+\sqrt{3}}\right)^{\chi_{12}(k)},
\end{equation}
where $\chi_{12}(l) = \left\{\begin{array}{lr}
        ~~\,1~\text{if} ~l\equiv \pm 1~ (\text{mod}~ 12)\\
        -1~\text{if} ~l\equiv \pm 5~ (\text{mod}~ 12)\\
        ~~~~~~~\,0~\text{otherwise}
        \end{array}\right\}$ is Dirichlet character modulo $12$.
        
$\mathbf{19}.$ If 
\begin{equation}
    \cosh\alpha_j+\cos\frac{\pi j}{4n}=\cosh\beta_k+\cos\frac{\pi k}{4m},\qquad (1\le j\le 4n,\ 1\le k\le 4m),
\end{equation}
then
\begin{equation}
    \prod_{j=1}^{4n}\left(\frac{\sqrt{2}\cosh m\alpha_j-1}{\sqrt{2}\cosh m\alpha_j+1}\right)^{\chi_{8}(j)}= \prod_{k=1}^{4m}\left(\frac{\sqrt{2}\cosh n\beta_k-1}{\sqrt{2}\cosh n\beta_k+1}\right)^{\chi_{8}(k)},
\end{equation}
where $\chi_{8}(l) = \left\{\begin{array}{lr}
        ~~\,1~\text{if} ~l\equiv \pm 1~ (\text{mod}~ 8)\\
        -1~\text{if} ~l\equiv \pm 3~ (\text{mod}~ 8)\\
        ~~~~~~~\,0~\text{otherwise}
        \end{array}\right\}$ is Dirichlet character modulo $8$.
        
Formulas in sections $17$ and $18$ are finite versions of propositions $27$ and $28$ in \cite{cais}.
        
$\mathbf{20}.$ If 
\begin{equation}
    \cosh2\alpha_j+\cos\frac{2\pi j}{3n}=\cosh2\beta_k+\cos\frac{2\pi k}{3m},\qquad (1\le j\le 3n,\ 1\le k\le 3m),
\end{equation}
then
\begin{equation}
    \prod_{j=1}^{3n}\left(1-\frac{3\sinh m\alpha_j}{\sinh 3m\alpha_j}\right)^{c(j)}=     \prod_{k=1}^{3m}\left(1-\frac{3\sinh n\beta_k}{\sinh 3n\beta_k}\right)^{c(k)},
\end{equation}
where $c(l)=2\cos\frac{2\pi l}{3}$.

Infinite version of this product is 
$$
\alpha\prod _{n=1}^\infty \left(1-\frac{3 \sinh\frac{\alpha n}{3}}{\sinh\alpha n}\right)^{c(n)}=
\beta\prod _{n=1}^\infty \left(1-\frac{3 \sinh\frac{\beta n}{3}}{\sinh\beta n}\right)^{c(n)}, \qquad\alpha\beta=\pi^2.
$$

{\it{Acknowledgements.}} The author of this paper wish to thank Dr. Lawrence Glasser for valuable correspondence and comments.

\end{document}